\pdfoutput=1
\documentclass[11pt]{article}
\usepackage[citecolor=blue]{hyperref}
\usepackage{amsmath,amsfonts,amssymb,amsthm}
\usepackage{mathrsfs}
\usepackage{verbatim}           
\usepackage{enumerate}
\usepackage[pdftex]{graphicx}
\usepackage[pdftex]{color}         
\usepackage{amscd}
\usepackage[footnotesize]{caption}

\numberwithin{equation}{section}
\def\yobs{{y_{\textrm{obs}}}}
\def\SNR{{\textrm{SNR}}}

\def\Z{{\mathbb Z}}
\def\R{{\mathbb R}}

\def\W{{\mathcal{W}}}

\def\L2{{L^2(\mathbb{R})}}

\def\Cov{{\mathrm{Cov}}}        
\def\E{{\mathbb{E}}}            

\def\|{{|\!|}}
\def\<{{\langle}}
\def\>{{\rangle}}

\title{A Multiresolution Ensemble Kalman Filter using Wavelet
Decomposition.}  
\author{ Kyle S. Hickmann$^{\thanks{ Applied Mathematics and Plasma
Physics, Los Alamos National Laboratory, hickmank@lanl.gov}}$,
Humberto C. Godinez$^{\thanks{ Applied Mathematics and Plasma Physics,
Los Alamos National Laboratory, hgodinez@lanl.gov}}$ }

\date{}
\begin{document}
\maketitle

\begin{abstract}
  We present a method of using classical wavelet based multiresolution
  analysis to separate scales in model and observations during data
  assimilation with the ensemble Kalman filter. In many applications, the
  underlying physics of a phenomena involve the interaction of features at
  multiple scales. Blending of observational and model error across scales can
  result in large forecast inaccuracies since large errors at one scale are
  interpreted as inexact data at all scales. Our method uses a transformation
  of the observation operator in order to separate the information from
  different scales of the observations. This naturally induces a
  transformation of the observation covariance and we put forward several
  algorithms to efficiently compute the transformed covariance. Another
  advantage of our multiresolution ensemble Kalman filter is that scales can
  be weighted independently to adjust each scale's effect on the forecast. To
  demonstrate feasibility we present applications to a one dimensional
  Kuramoto-Sivashinsky (K--S) model with scale dependent observation noise and
  an application involving the forecasting of solar photospheric flux. The
  latter example demonstrates the multiresolution ensemble Kalman filter's
  ability to account for scale dependent model error. Modeling of photospheric
  magnetic flux transport is accomplished by the Air Force Data Assimilative
  Photospheric Transport (ADAPT) model.
\end{abstract}

\footnotetext{{\bf 2010 MSC}:}

\medskip


\medskip


\section{Introduction}
\label{sec:introduction}

Combining large-scale physics simulations with data to generate informed
forecasts, with quantified uncertainty, is a common task in modern science
\cite{Daley1991,evensen_data_2009,kalnay_atmospheric_2003}. A prevalent method
to accomplish this task is the ensemble Kalman filter (EnKF)
\cite{evensen1994sequential} which can provide a forecast of the mean behavior
of the system along with confidence intervals. The EnKF assumes that both the
model and observations have a Gaussian probability distribution.
Additionally, it is frequently assumed that model and observation covariance
matrices are diagonal or block diagonal, implying that errors associated with
different state variables are uncorrelated or weakly correlated. These types
of assumptions are even more prevalent when variables or observations are far
apart in space and/or time \cite{HoutekamerMitchell1998}. However, if
observation or model errors are scale dependent, correlations of errors
between variables may not be known \emph{a priori} and discarding some of
these correlations artificially can cause ensemble collapse resulting in large
forecast errors.

The problem of scale dependent observation and model error in forecasting can
be seen in atmospheric data assimilation where models that do a good job of
forecasting large scale phenomenon are coupled with models of small scale
turbulent effects
\cite{casati2004new,chou1993multiresolution,harris2001multiscale,palmer2001nonlinear}.
Here we propose to address the scale dependence problem in the EnKF by using a
multiresolution analysis (MRA).  Wavelet based MRA has been used to analyze
the statistical properties of weather models and ocean models in the past
\cite{beezley2011wavelet,buehner2012evaluation,buehner2007spectral,chou1993multiresolution,deckmyn2005wavelet,kasanicky2014spectral}.
However, this has usually been done retrospectively. Here we propose to use
the results of MRA scale separation during the generation of a data
assimilative forecast by directly coupling the MRA with the EnKF.

Wavelet analysis has been applied to the EnKF in the past
\cite{beezley2011wavelet,buehner2007spectral,deckmyn2005wavelet}. Previously
MRA was used to offer a more accurate approximation to the ensemble
covariance. Past work first transformed the ensemble to the wavelet domain
where a diagonal wavelet covariance was estimated, before assimilation was
performed, the wavelet covariance was then transformed back to the original
model domain. Approximating the covariance in this way offers a more accurate
estimation of the ensemble covariance due to a regularization effect that the
wavelet transform naturally provides. However, the improvement in this
approximation relies on the assumption that the wavelet transform
approximately diagonalizes the ensemble covariance, which is not always the
case \cite{dijkerman1994wavelet,masry1998covariance}.

Here, we apply the wavelet transform to the observation operator directly. The
transformation to the wavelet domain is computed only once during assimilation
as a preprocessing step and the inverse transformation is never
computed. Using the wavelet transform to modify the observation operator,
instead of modifying the ensemble covariance, has the effect of offering a
computationally efficient, scale dependent, extension of the EnKF.  We also
show that an iterative application of the EnKF with a scale dependent
observation operator allows for propagation of information between
scales. Thereby eliminating the need for assumptions about independence of
scales. Another advantage of our method is that, once the transformation to
the wavelet domain has been computed, it is natural to use an ensemble
inflation coefficient to assign trust to the observations and model error
based on \emph{a priori} knowledge about the accuracy of observations and
model at each scale.

We demonstrate our methods on two different models. First, we apply the
multiresolution EnKF (MrEnKF) to the Kuramoto-Sivashinsky (K--S) equation. The
K--S equation is a 1D nonlinear partial differential evolution equation which
possesses multi-scale dynamics. The K--S equation is used here to demonstrate
the advantage of the MrEnKF in a data assimilation experiment in which we
assign varying degrees of observational noise to distinct wavelet scales. It
is assumed that the large scales are observed more accurately than the finer
scales. Since the large scales contain more information about the unstable low
Fourier modes, this observation experiment demonstrates the advantages of the
MrEnKF over the scale independent EnKF when dealing with multi-scale models.

The second example we present involves forecasting the magnetic flux
transported across the solar photosphere.  Our application is the Air Force
Data Assimilative Photospheric Transport (ADAPT) model of photospheric flux
propagation
\cite{arge2013modeling,arge2010air,arge2011improving,hickmann2015data}, a
collaborative modeling and forecasting effort between Los Alamos National
Laboratory and the Air Force Research Laboratory in Kirkland ABQ.  The solar
photosphere application highlights the challenges encountered in realistic
modeling and forecasting efforts within the science community.  Many scientist
are interested in the tracking of emergent coherent regions of magnetic
flux. These large clumps of magnetic flux are known as \emph{active regions}
and are primary drivers of large space weather events such as Coronal Mass
Ejections (CMEs)
\cite{antiochos1999model,falconer2002correlation,glover2000onset,munoz2010double}.
When implementing a standard EnKF, the active regions tend to diffuse and lose
structure after only one assimilation cycle. By the end of the assimilation
window, it is difficult to maintain a coherent active region structure with
the EnKF. On the other hand, we show that the MrEnKF performs much better in
regards to maintaining a coherent active region, since the MrEnKF assigns
greater confidence to observations characteristics at the scale of active
regions. Once the structure is preserved for a newly emerging active region,
successive observations of the active region allow for increasing definition
in ADAPT's data assimilation mechanism.

In section 2 we quickly review the classical wavelet multiresolution analysis,
set up notation, and point the interested reader to references on wavelet
analysis. Section 3 then gives a derivation of our multiresolution ensemble
Kalman Filter scheme. In section 4 we discuss ways to approximate the change
to the observation covariance when using the wavelet transformation in a
computationally efficient manner. Section 5 details the role of ensemble
inflation in the multiresolution EnKF. Our examples using the
Kuramoto-Sivashinsky equation and the ADAPT forecasting model are detailed in
section 6. We conclude with a discussion of hopeful applications of the MrEnKF
as well as future improvements.

\section{Wavelet Decomposition}
\label{sec:wavelet}

Wavelet analysis has been used in a wide spectrum of applications where
fidelity of information varies by location and frequency simultaneously or
where one seeks to isolate a particular signal in both location and
frequency. Wavelet analysis has its roots in Fourier analysis, where one
decomposes a signal with respect to frequency.  However, in Fourier analysis
the frequency information is not localized in the original domain, which can
be undesirable when performing time-series analysis or image processing. For
this reason wavelet analysis, a localized frequency decomposition, was
developed.  The goal of this section is to establish notation that will be
used throughout the paper and briefly revisit some concepts of wavelet
decomposition. The wavelet analysis used in our work is based on the
multiresolution decomposition of Mallat
\cite{mallat1989multiresolution,mallat1989theory}. For further details
regarding wavelets we refer the reader to the work by Daubechies
\cite{daubechies1992ten}

Suppose we are given a discretized signal $f \in \R^n$. We can perform an
$N$-level wavelet transform \cite{mallat1989multiresolution} of the observed
signal, $w^f = \W_N f$ and group the coefficients by level,
\begin{equation*}
  w^f = \left[ (w^f_{N+1})^T, (w^f_N)^T, (w^f_{N-1})^T, \cdots, (w^f_1)^T \right]^T. 
\end{equation*}
Here $w^f_{N+1}$ represents the coarsest coefficients and each successive
$w^f_i$, $i = N, N-1, \dots, 1$, represents increasingly fine scale
coefficients. We define the projections of the wavelet coefficients onto each scale by
\begin{equation}\label{wav_projection}
  P_i w^f = w^f_i, \qquad i = N+1, N, N-1, \dots, 1.
\end{equation}
In general the coarser wavelet coefficients capture larger scale behavior of
the signal, each successive wavelet level captures finer scale variations
\cite{daubechies1992ten,mallat1989multiresolution,mallat1989theory}.

For purposes of separating scales within the ensemble Kalman filter, it is
necessary to compute the effect of the multi-scale decomposition on the
covariance of both the model and the observation. If we assume that a signal
has a Gaussian distribution, $f \sim N( \mu_f, C_f)$, then the wavelet
coefficients are Gaussian distributed as well, $w^f_i \sim N( \mu^f_i, C^f_i)$
for $i = N+1, \dots, 1$. The mean and covariance for each level of wavelet
coefficients is given by
\begin{equation*}
  \mu^f_i = P_i \W_N \mu_f,\,\, C^f_i = P_i \W_N C_f (P_i \W_N)^T.
\end{equation*}
This result relies only on the linearity of the wavelet transform and
projections. Here we have given formulas only for the covariance matrices for
each level of the $N$-level wavelet transform, ignoring covariance terms
between levels of the transform. It is worth noting that our proposed MrEnKF
allows some interaction between scales though these interactions are not
treated explicitly.

\section{Multiresolution ENKF}
\label{sec:mrenkf}

Using the MRA decomposition introduced in the previous section we put forward
a method of including scale dependent model and observation error information
in the EnKF assimilation scheme. The method we propose is iterative, over
scales, which allows for a limited increase in computational complexity over
the standard EnKF methods. Moreover, the method we use to include the wavelet
decomposition modifies the observation operator of the assimilation problem
and therefore is agnostic toward the exact EnKF implementation,
e.g. stochastic, square root, transform etc.

The common setup for the EnKF is as follows
\cite{Daley1991,evensen_data_2009,kalnay_atmospheric_2003}, an observation of
a dynamical system is made with assumed Gaussian observational error, $\yobs
\sim N(\mu_{\textrm{obs}}, R)$. The forecast (background) state, denoted
$x_b$, is assumed to be generated by a mathematical model of the system. Due
to model error the forecast state is assumed to be Gaussian distributed, $x_b
\sim N(\mu_b, C_b)$. The forecast state is related to the observation through
a, possibly nonlinear, observation operator $H(\cdot)$ and we write the
forecast observation as $H(x_b)$. The model is assumed to be unbiased so that
$\mu_{\textrm{obs}} = \E[H(x_b)]$. Under these assumptions the pair $(x_b,
\yobs)$ is jointly Gaussian distributed,
\begin{align}
  \left[\begin{array}{c}
      x_b \\ \yobs
  \end{array}\right] \sim N\left( 
\left[\begin{array}{c}
\mu_b \\ \E[H(x_b)]
\end{array}\right],
\left[\begin{array}{cc}
  C_b & \Cov(x_b, H(x_b)) \\
  \Cov(H(x_b), x_b) & \Cov(H(x_b)) + R
\end{array}\right]
\right).
\end{align}
We refer to the forecast state conditioned on the observations as the
\emph{analysis} state, $x_a = x_b | \yobs$. Given the preceding assumptions
the analysis has Gaussian distribution and we write $x_a \sim N(\mu_a, C_a)$
with
\begin{equation}\label{enkf}
\begin{aligned}
  \mu_a &= \mu_b + \Cov(x_b, H(x_b))(\Cov(H(x_b)) + R)^{-1} (\yobs - \E[H(x_b)]) \\
  C_a &= C_b - \Cov(x_b, H(x_b))(\Cov(H(x_b)) + R)^{-1} \Cov(H(x_b), x_b).
\end{aligned}
\end{equation}
The analysis state of the system represents the distribution of likely system
states given our most recent observations. In order to approximate draws of
$x_a$ the EnKF first generates a forecast \emph{ensemble} of states
$x^{\alpha}_b$ for $\alpha = 1, 2, \dots, M$, where $M$ is the size of the
ensemble. The forecast ensemble is used to compute sample expectations and
covariances to replace the terms $\mu_b$, $C_b$, $\E[H(x_b)]$, $\Cov(x_b,
H(x_b))$, and $\Cov(H(x_b))$ in (\ref{enkf}). Once (\ref{enkf}) has been
approximated using the forecast ensemble there are a myriad of methods to
generate an analysis ensemble $x^{\alpha}_a$, $\alpha = 1, 2, \dots, M$,
approximating draws from $N(\mu_a, C_a)$. These methods include stochastic
EnKF variations \cite{evensen1994sequential,evensen1996assimilation}, ensemble
transform Kalman Filters (ETKF) \cite{hunt2007efficient}, ensemble adjustment
Kalman Filters \cite{anderson2001ensemble}, and ensemble square root Kalman
filters \cite{tippett2003ensemble,whitaker2002ensemble}. We refer to all of
these methods loosely as ensemble Kalman filters. In our applications we use
an ensemble transform Kalman filter \cite{hunt2007efficient}. Once the
analysis ensemble is generated a new forecast is computed by propagating the
analysis ensemble through the mathematical model until the next observation
time and the whole process is iterated.

In order to include multi-scale information into the EnKF, we apply the
wavelet decomposition to the observation operator $H(\cdot)$. In practice this
means that the wavelet transform is applied to the data and ensemble of state
observations as a preprocessing step. Decomposing $H(\cdot)$ yields observed
wavelet coefficients $w_{\textrm{obs}} = \W_N \yobs$ and the unbiased
observation assumption implies $w_{\textrm{obs}} = \E[\W_N H(x_b)]$. Scales
are separated by applying the wavelet projections (\ref{wav_projection}),
yielding $N+1$ sets of scale dependent observations, 
\begin{equation}
  y^i_{\textrm{obs}} = P_i w_{\textrm{obs}}, \qquad i = N+1, N, \dots, 1
\end{equation}
with observation errors distributed as
\begin{equation}
  y^i_{\textrm{obs}} \sim N(P_i \W_N \mu_{\textrm{obs}},  P_i \W_N R (P_i \W_N)^T).
\end{equation}
For convenience we define the scale dependent observation operator as 
\begin{equation}
H_i(x_b) = P_i \W_N H(x_b), \qquad i = N+1, N, \dots, 1
\end{equation}
and note that since the model is assumed unbiased
\begin{equation}
  \E[y^i_{\textrm{obs}}] = \E[ H_i(x_b) ], \qquad i = N+1, N, \dots, 1.
\end{equation}
With the above notation we may express the unbiased assumption at each scale
as
\begin{equation}
  y^i_{\textrm{obs}} = H_i(x_b) + \epsilon_i, \qquad \epsilon_i \sim N(0, R_i)
\end{equation}
where $R_i = P_i \W_N R (P_i \W_N)^T$ for $i = N+1, N, \dots, 1$. 

Instead of conditioning $x_b$ on all observation scales simultaneously we
iteratively condition the forecast on the observations from one scale at a
time. Largest scales are assimilated first followed by the assimilation of
finer scales in the observations. Of course any other ordering is possible and
perhaps the most accurately observed scales should be assimilated first.

We use the notation $x_{N+2} = x_b$, $\mu_{N+2} = \mu_b$, and $C_{N+2} =
C_b$. Our iterative multiresolution EnKF (MrEnKF) is then defined by the
series of conditioned model states $x_i = x_{i+1} | y^i_{\textrm{obs}} \sim
N(\mu_i, C_i)$ for $i = N+1, N, \dots, 1$ with mean and covariance at each
scale given by
\begin{equation}\label{mrenkf}
\begin{aligned}
  \mu_i &= \mu_{i+1} + \Cov(x_{i+1}, H_{i}(x_{i+1})) [ \Cov(H_i (x_{i+1})) + R_i]^{-1} (y^i_{\textrm{obs}} - \E[H_i(x_{i+1})]) \\
C_i &= C_{i+1} - \Cov(x_{i+1}, H_i (x_{i+1})) [ \Cov(H_i (x_{i+1})) + R_i]^{-1} \Cov(H_i (x_{i+1}), x_{i+1}).
\end{aligned}
\end{equation}
The final analysis, with all scales assimilated, then is distributed as $x_a =
x_{1} = x_{2} | y^1_{\textrm{obs}} \sim N(\mu_1, C_1)$.

In practice each of the Gaussian distributions is approximated by an ensemble
of model states. Any standard EnKF type algorithm can be used to form an
analysis ensemble approximating a draw from the Gaussian conditioned by each
successive scale. In our examples we have used the ensemble transform Kalman
filter (ETKF) \cite{hunt2007efficient} to form the intermediate analysis
ensembles but any EnKF variation could be applied. It is important to note two
properties of our MrEnKF. First, we have made the important assumption that
the \emph{conditioning on the larger scale does not effect the bias at the
  finer scale}. Second, because the conditioning on scales is performed
iteratively we did not need to make any assumptions about independence of
scales and scale/scale covariance information is allowed to propagate through
to the analysis.

\section{Observation Covariance}
\label{sec:obscov}

For large problems the observation covariance transformation may be very
expensive to compute at each level. Computing $R_i = P_i \W_N R (P_i \W_N)^T$
for $i = N+1, N, \dots, 1$ requires calculating a wavelet transform for each
dimension of the model forecast $x_b$. In cases where the data being
assimilated is the discretization of a two or three dimensional field
computing and storing this covariance matrix is computationally
prohibitive. This is further complicated by the fact that the wavelet
transform is not usually stored as a matrix, so computing $(P_i \W_N)^T$ is
not straight forward. Some of these problems can be simplified if we know the
symmetric square root decomposition of the covariance matrix, $R = S S^T$. If
the square root is available then at least the question of transposing the
wavelet transform is averted and $R_i = (P_i \W_N S) (P_i \W_N S)^T$.

If the original covariance $R$ is diagonal or very dominated by the diagonal
terms we may be safe in making the assumption that $R_i$ is a diagonal matrix
with a constant on the diagonal determined by some overall measure of the
noise level in the observations at each scale. This can be accomplished by
setting
\begin{equation}
R_i = \lambda_i \,\, \sigma^2_{\textrm{max}} (R), \qquad i = N+1, N, \dots, 1.
\end{equation}
Here $\sigma_{\textrm{max}} (R)$ is the largest singular value of $R$ and $0 <
\lambda_i$ is a scaling parameter to adjust the confidence given to
observations at each scale. The latter approach has worked well in our
examples (Section \ref{sec:examples}) but requires tuning of the scaling
parameter and represents a drastic assumption about the observation error
within each scale.

A more rigorous way of approximating the covariance for each observation scale
relies on the sampling of observations and approximating the covariance at
differing scales in the spirit of the original EnKF proposed by Evensen
\cite{evensen1994sequential,evensen1996assimilation}. This method is accurate
if one is willing to sample the observational noise determined by $R$. The
method starts by generating $M$ samples of $\epsilon \sim N(0, R)$,
$\epsilon_j, j = 1, 2, \dots, M$. To approximate $R_i$ for a given wavelet
scale each $\epsilon_j$ is transformed to give $\epsilon^i_j = P_i \W_N
\epsilon_j$. The transformed noise samples are then used to form a noise
ensemble matrix
\begin{equation}
E_i = [\epsilon^i_1 | \epsilon^i_2 | \cdots | \epsilon^i_M].
\end{equation}
The covariance is then approximated by 
\begin{equation}
  R_i \approx \frac{1}{M-1} E_i E^T_i.
\end{equation}
The disadvantage of this approach is that we may require a large number of
noise samples to accurately approximate $R_i$. However, if covariance
inflation is to be used at each scale, a very accurate approximation may not
be necessary. Moreover, at least the larger scale components will have a
significantly lower dimension than the original forecast, and therefore will
allow an accurately approximated covariance with far fewer samples than would
be necessary to approximate the full covariance.

\section{Ensemble Inflation}
\label{sec:inflation}

Ensemble inflation has been shown to be beneficial in preventing ensemble
collapse and divergence when using ensemble data assimilation schemes
\cite{HoutekamerMitchell1998,wang2003comparison}. Moreover, in
\cite{kelly2014well} it was shown that in order to have both stability and
accuracy in an EnKF scheme inflation was necessary. The MrEnKF proposed makes
it straightforward to apply a scale dependent inflation, making for a very
robust/tunable filter. The analysis mean and covariance at each scale given in
equation (\ref{mrenkf}) can be replaced with {\small
\begin{equation}\label{inflation_mrenkf}
\begin{aligned}
  \mu_i &= \mu_{i+1} + \Cov(x_{i+1}, H_i (x_{i+1})) \left[ \Cov(H_i (x_{i+1})) + \frac{1}{\rho_i} R_i \right]^{-1} (y^i_{\textrm{obs}} - \E[H_i (x_{i+1})]) \\
C_i &= C_{i+1} - \Cov(x_{i+1}, H_i (x_{i+1})) \left[ \frac{1}{\rho_i} \Cov(H_i (x_{i+1})) + \frac{1}{\rho^2_i} R_i \right]^{-1} \Cov(H_i (x_{i+1}), x_{i+1})
\end{aligned}
\end{equation}
} for $i = N+1, N, \dots, 1$. This implies that a vector of scale dependent
inflation coefficients must be chosen,
\begin{equation}
\rho = (\rho_{N+1}, \rho_{N}, \dots, \rho_{1})^T.
\end{equation}
At each scale the coefficient $0 < \rho_i$ quantifies the amount of confidence
given to either the model or the observation during assimilation.

A scale dependent inflation allows the user to control the confidence in the
model or observation at each scale separately. Therefore, if a set of
observations is known to be a very accurate measure of one scale the inflation
coefficient for this scale can be increased while the others are left
unchanged. Allowing this level of tuning can be advantageous if there is
detailed information available about the scale dependence of observation and
model errors.

\section{Example Applications}
\label{sec:examples}

We give two examples where the use of an MrEnKF scheme is beneficial. Scale
dependence can become important if the observation error varies greatly with
scale or if the model error varies greatly with scale. In our first example we
apply the MrEnKF to a chaotic nonlinear PDE in one dimension under the
assumption of scale dependent observational noise. The second example
demonstrates the MrEnKF on a problem in solar weather forecasting in which the
model itself has scale dependent error. In both cases we demonstrate that
separation of scales during the assimilation can significantly improve the
ensemble's ability to track the observed data while reliably representing the
error in the forecast.

\subsection{Kuramoto-Sivashinsky Equation}
\label{sec:KS}

The Kuramoto-Sivashinsky equation was named for its derivation in modeling
hydrodynamic stability of laminar flame fronts \cite{sivashinsky1977nonlinear}
and as a phase equation for the complex Ginzburg-Landau equation
\cite{kuramoto1976persistent}. It was first derived as a model of nonlinear
saturation of drift waves associated with the oscillation of plasma particles
trapped in magnetic wells \cite{laquey1975nonlinear}. Applications of the K--S
equation include modeling of the dynamics of self-focusing lasers
\cite{munkel1996intermittency}, instabilities in thin films
\cite{babchin1983nonlinear}, and the flow of a viscous fluid down a vertical
plane \cite{sivashinsky1980irregular}. Extensive numerical investigations of
the chaotic dynamics of the K--S equation have been carried out
\cite{drotar1999numerical,hyman1986kuramoto,hyman1986order,khellat2014kuramoto,smyrlis1996computational}. Furthermore,
the K--S equation has been a source of many results related to dynamics of
chaotic systems
\cite{collet1993global,elgin1996stability,nicolaenko1985some}. In regards to
assimilation and control the K--S equation has classically represented a
challenging problem to test methods of control and assimilation for chaotic
dynamical systems
\cite{desertion2004improved,dubljevic2010model,el2008actuator,hu2001robust,jardak2010comparison}. For
an overview of the theory of existence and uniqueness of the K--S equation the
interested reader is pointed to \cite{robinson2001infinite,temam2012infinite}.

In its simplest form the Kuramoto-Sivashinsky equation in one dimension is expressed as 
\begin{equation}\label{KSeqtn}
\begin{aligned}
  u_t + u_{xx} &+ u_{xxxx} + u u_x = 0 \quad \textrm{ on } [-\pi L, \pi L] \times [0, T] \\
  u(x,0) &= u_0(x),\, u(-\pi L, t) = u(\pi L, t) \quad \textrm{ for } t \ge
  0,
\end{aligned}
\end{equation}
where $L$ is a bifurcation parameter which controls the behavior of solutions,
i.e. stable, periodic, chaotic, etc. Equation (\ref{KSeqtn}) is diagonalized
by the Fourier transform to get the system of ODEs
\begin{equation}
\begin{aligned}
     u(x,t) &= \sum_n u_n(t) \exp \left( \frac{inx}{L} \right), \,\,\, u_0(x) = \sum_n u_n(0) \exp \left( \frac{inx}{L} \right) \\
     \frac{du_n}{dt} &= \left( \frac{n}{L} \right)^2 \left(1 - \left( \frac{n}{L} \right)^2 \right) u_n - \frac{in}{2L} \sum_{j \in \Z} u_j(t) u_{n-j}(t).
\end{aligned}
\end{equation}
This diagonalization shows that the first $0 \le n < L$ Fourier modes are
unstable about $u_n(t) = 0$ while the higher Fourier modes are stable
\cite{kassam2005fourth,robinson2001infinite,temam2012infinite}. The nonlinear
term then allows mixing between the low and high Fourier modes which allows
for stable solutions as some of the energy is transferred from the low to the
high modes and then dissipated
\cite{robinson2001infinite,temam2012infinite}. This property of the K--S
equation makes it ideal for testing a scale dependent EnKF since we can assume
that the unstable low frequencies, large scales, are observed with higher
accuracy than the high frequencies, small scales, and investigate the effect
of propagating scale dependent information through the EnKF.

In our data assimilation experiments with the K--S equation we assume
$L = 22$ which is well into the regime of chaotic solutions
\cite{jardak2010comparison,kassam2005fourth}. Our experiment's initial
condition's are chosen as in the work of
\cite{jardak2010comparison,kassam2005fourth},
\begin{equation}\label{KS_IC}
  u_0(x) = \cos \left( \frac{x}{L} \right) \left(1 + \sin \left(\frac{x}{L} \right)\right).
\end{equation}
Solutions to the K--S equation are simulated using a stable fourth
order Runge-Kutta scheme with exponential time differencing
\cite{jardak2010comparison,kassam2005fourth}. The spatial domain is
discretized using $512$ equally spaced points on $-\pi L < x_j \le \pi
L$, the temporal domain is discretized using a step length of $\Delta
t = 0.5$. We assume that solutions of (\ref{KSeqtn}) are observed,
until time $T = 300$, every $20^{th}$ time-step $t_n = 20 n \Delta t$.
Observations of the K--S solution will be denoted by $H(u(t_n)) =
(u(x_1,t_n), u(x_2, t_n), \dots, u(x_{512}, t_n))^T$.

We set up a \emph{twin} experiment with the K--S equation to compare the
performance of the EnKF and MrEnKF. The twin experiment consists of simulating
observations from a reference solution of the K--S equation having initial
conditions given by (\ref{KS_IC}). Observations are generated from this
reference solution at each of the time points $t_n$ and Gaussian noise, as
specified in the next paragraph, is added to the observations. For both the
EnKF and MrEnKF ensembles are initialized by adding Gaussian white noise with
standard deviation $\sigma = 0.8$ to the initial condition (\ref{KS_IC}). The
ensemble size in all experiments was taken to be $N = 50$. The ensemble
members are propagated forward according to the K--S equation to time $t_1$
and assimilation is performed using the synthetic observations. This is
repeated at each $t_n$ and the results of the assimilation are compared.

Scale dependent observation error is modeled by taking a level-$4$ wavelet
transform of the solution at each $t_n$ using a Daubechies 'db9' wavelet with
$9$ vanishing moments \cite{daubechies1992ten}.  Gaussian white noise is then
added to each level of wavelet coefficients with standard deviation dependent
on the transform level of the coefficients. The standard deviation of the
additive noise at each wavelet level was chosen so that the signal-to-noise
ratio ($\SNR$) was smaller for the fine scale coefficients, thus keeping more
of the large scale information in the observations. Specifically, the standard
deviations for the Gaussian noise added to the four levels of wavelet
coefficients were taken to be $\sigma_{5} = 0.75$, $\sigma_{4} = 0.75$,
$\sigma_{3} = 1.65$, $\sigma_{2} = 1.0$, $\sigma_{1} = 0.0008$. This choice of
standard deviation gave an average $\SNR$ of $\SNR_5 \approx 18.22$, $\SNR_4
\approx 15.58$, $\SNR_3 \approx 2.04$, $\SNR_2 \approx 1.16$, and $\SNR_1
\approx 1.17$ at the respective scales. Average SNR values were calculated by
applying the formula
\begin{equation}
  \SNR_i = \frac{\max (w_i) - \min
  (w_i)}{\sigma_{\mathrm{noise}}}, \qquad i = N+1, N, \dots, 1
\end{equation}
to each scale of the wavelet coefficients and then averaging each scale's
$\SNR$ over all observation times. Here $\max(w_i)$ is the maximum wavelet
coefficient at level $i$, similarly $\min(w_i)$ is the minimum wavelet
coefficient at level $i$. The actual deviations of the simulated observations,
with this noise structure, from the true solution were Gaussian and had a
standard deviation of $\sigma_{\mathrm{obs}} = 0.8$.

\begin{figure}[h]
  \begin{center}
      \includegraphics[scale=0.6]{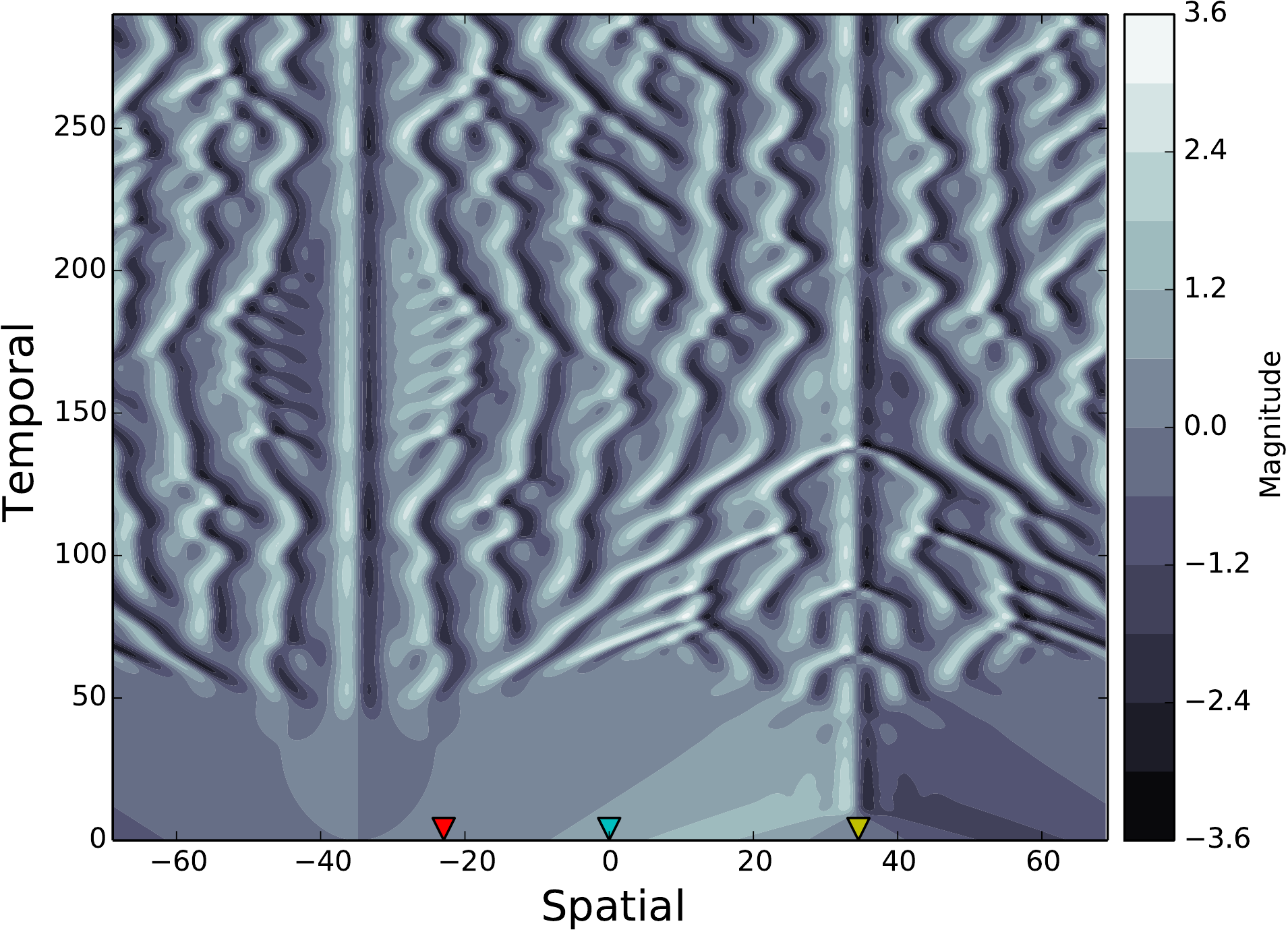}
  \end{center}
  \caption{Our reference solution for the K--S equation is shown as a contour
    plot. The red, cyan, and yellow markers indicate points where our ensemble
    forecasts are presented in Figure
    \ref{fig:etkfVSmretkf_pntwise}. Observations were generated by adding
    scale dependent Gaussian noise to this solution.}
  \label{fig:KSsoltn_markers}
\end{figure}

\begin{figure}[h]
  \begin{center}
    $\begin{array}{cc}
      \includegraphics[scale=0.37]{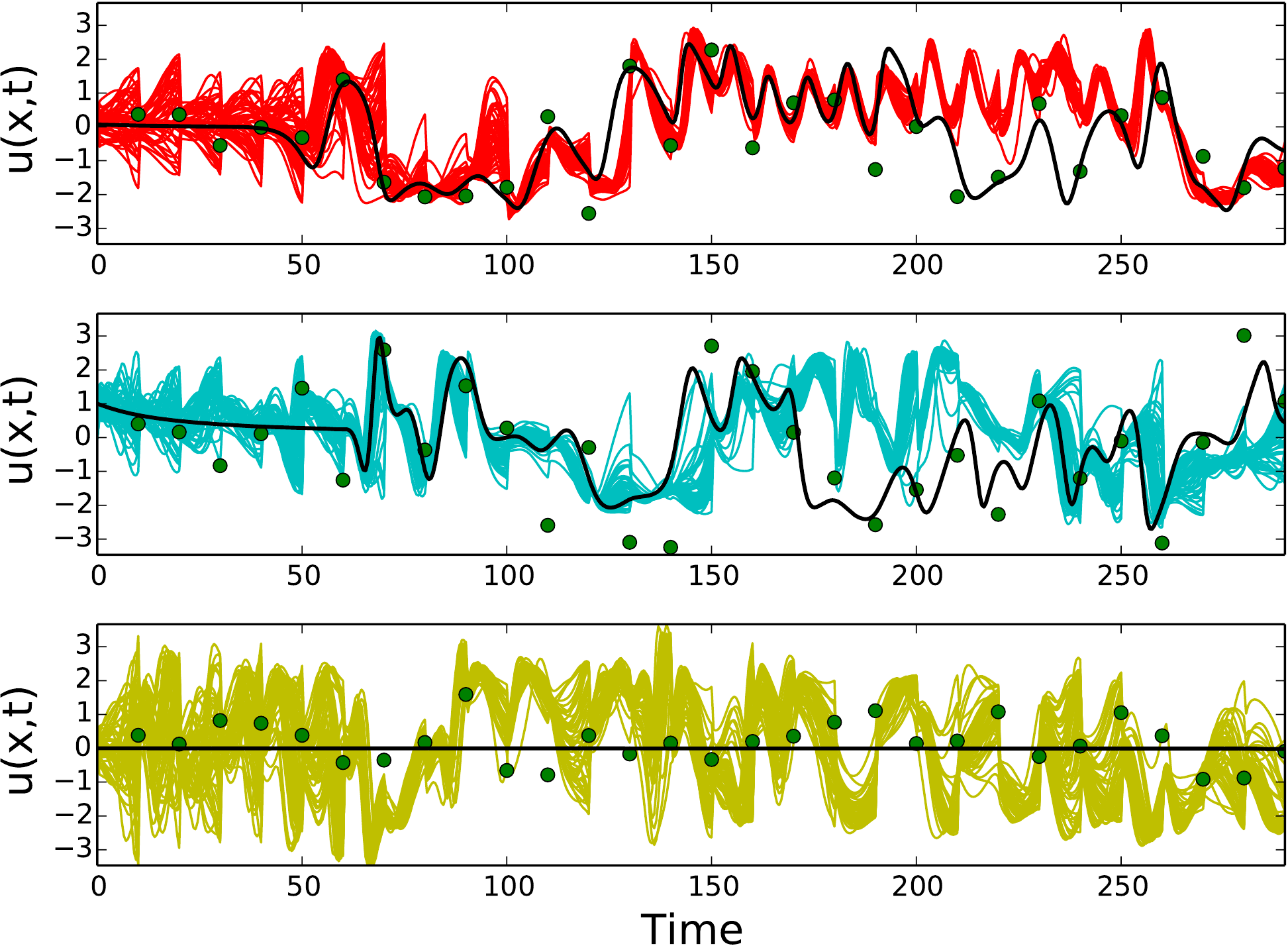} &
      \includegraphics[scale=0.37]{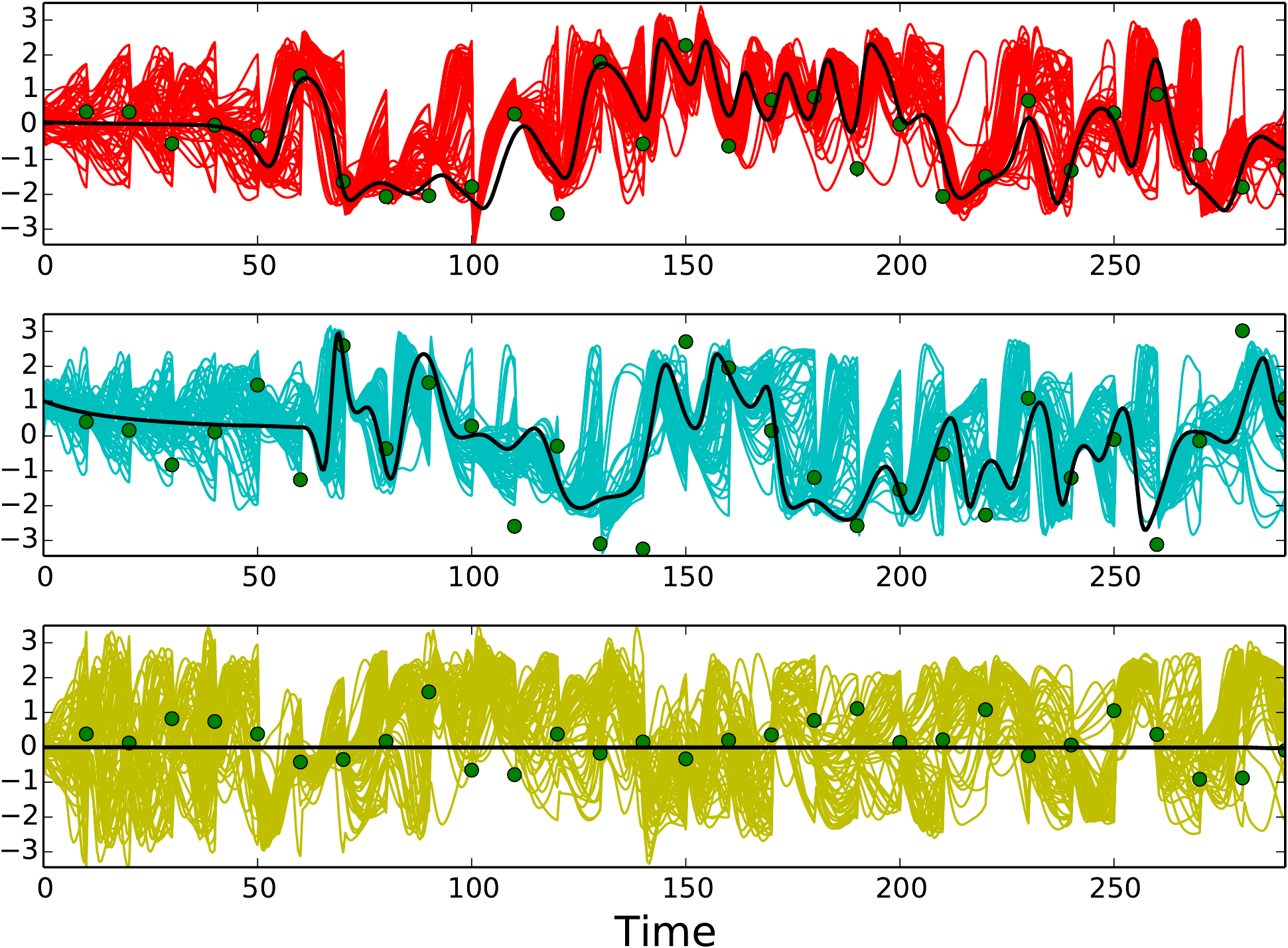}
      \end{array}$
  \end{center}
  \caption{Tracking of a solution of the K--S equation using the EnKF (LEFT)
    and MrEnKF (RIGHT). Ensemble colors correspond to the marker colors in
    Figure \ref{fig:KSsoltn_markers}. The true solution values at each point
    are shown in black while the solution observations are shown with green
    dots. Scale dependent noise has been added to the observations of this
    solution as described. We see that the ensemble, under EnKF, has
    significant periods of divergence from the true solution. More problematic
    is the fact that the standard deviation of the ensemble is small compared
    to the observational noise indicating a great deal of confidence in the
    assimilation. The MrEnKF on the other hand is shown to track the solution
    in the sense that usually the ensemble envelops the true solution. The
    standard deviation of the ensemble is kept large by taking into account
    scale dependent observation error. Due to the large ensemble spread the
    $L^2$ discrepancy, Figure \ref{fig:etkfVSmretkf_L2discrepancy}, for the
    MrEnKF is only incrementally better than for the EnKF. However, the
    standard deviation is more in line with the true observation error and the
    forecast spread reflects the variance in the true solution more accurately
    as demonstrated by the rank histograms in Figure
    \ref{fig:etkfVSmretkf_rankhistogram}}
  \label{fig:etkfVSmretkf_pntwise}
\end{figure}

The EnKF and the MrEnKF were applied to the Kuramoto-Sivashinsky assimilation
problem with scale independent noise. In this case the EnKF and the MrEnKF
both resulted in ensembles that tracked the solution very well for Gaussian
observation noise with a standard deviation of $\sigma_{\mathrm{obs}} =
0.8$. The EnKF and MrEnKF were then compared using the scale dependent
observation noise structure described above. To visualize the ensemble's
ability to track the true solution of the K--S equation we present the
ensemble tracking at three distinct points within the domain $[-\pi L, \pi
L]$. Figure \ref{fig:KSsoltn_markers}
shows the true solution of the K--S equation used to generate observations
along with three markers indicating the spatial points where ensemble tracking
is illustrated in Figure \ref{fig:etkfVSmretkf_pntwise}.
The K--S solution with initial condition (\ref{KS_IC}) has two stationary
nodes at $x = \pm \frac{L \pi}{2}$. We illustrate the ensemble tracking near
the first stationary node $\hat{x}_1 = -7.3 \pi$, away from both stationary
nodes $\hat{x}_2 = 0.0$, and on the second stationary node $\hat{x}_3 = 11.0
\pi$. At $\hat{x}_3$ both the EnKF and the MrEnKF have difficulties tracking
the solution. This is due to the fact that locally the K--S solution is not
stationary and therefore the incorporation of this local information in the
EnKF and the MrEnKF tends to pull the ensemble away from zero. However, the
MrEnKF still envelops the stationary solution much better, though with an
admittedly high ensemble standard deviation. At the points $\hat{x}_1$ and
$\hat{x}_2$ the EnKF reduces the ensemble spread very quickly and then can be
forced off of the true solution randomly (left plots in Figure
\ref{fig:etkfVSmretkf_pntwise}). The MrEnKF does not have this behavior (right
plots in Figure \ref{fig:etkfVSmretkf_pntwise}), though it must pay the cost
of maintaining a larger ensemble spread.  The MrEnKF ensemble almost always
envelops the true solution and therefore offers a significant advantage. The
advantages of the MrEnKF become clearer if we compare the rank histograms
\cite{hamill2001interpretation} for the EnKF and MrEnKF.

\begin{figure}[h]
  \begin{center}
    $\begin{array}{cc}
      \includegraphics[scale=0.38]{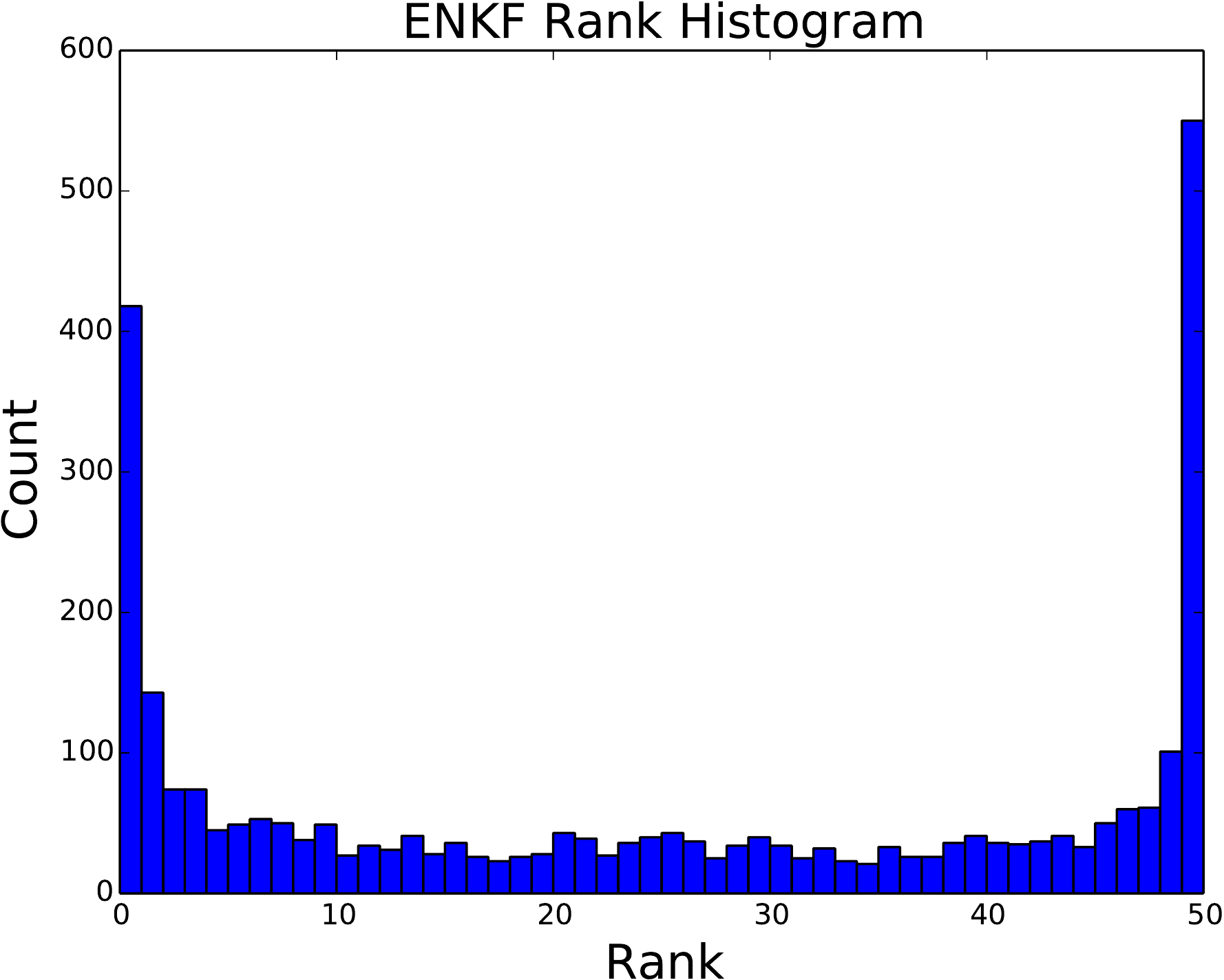} &
      \includegraphics[scale=0.38]{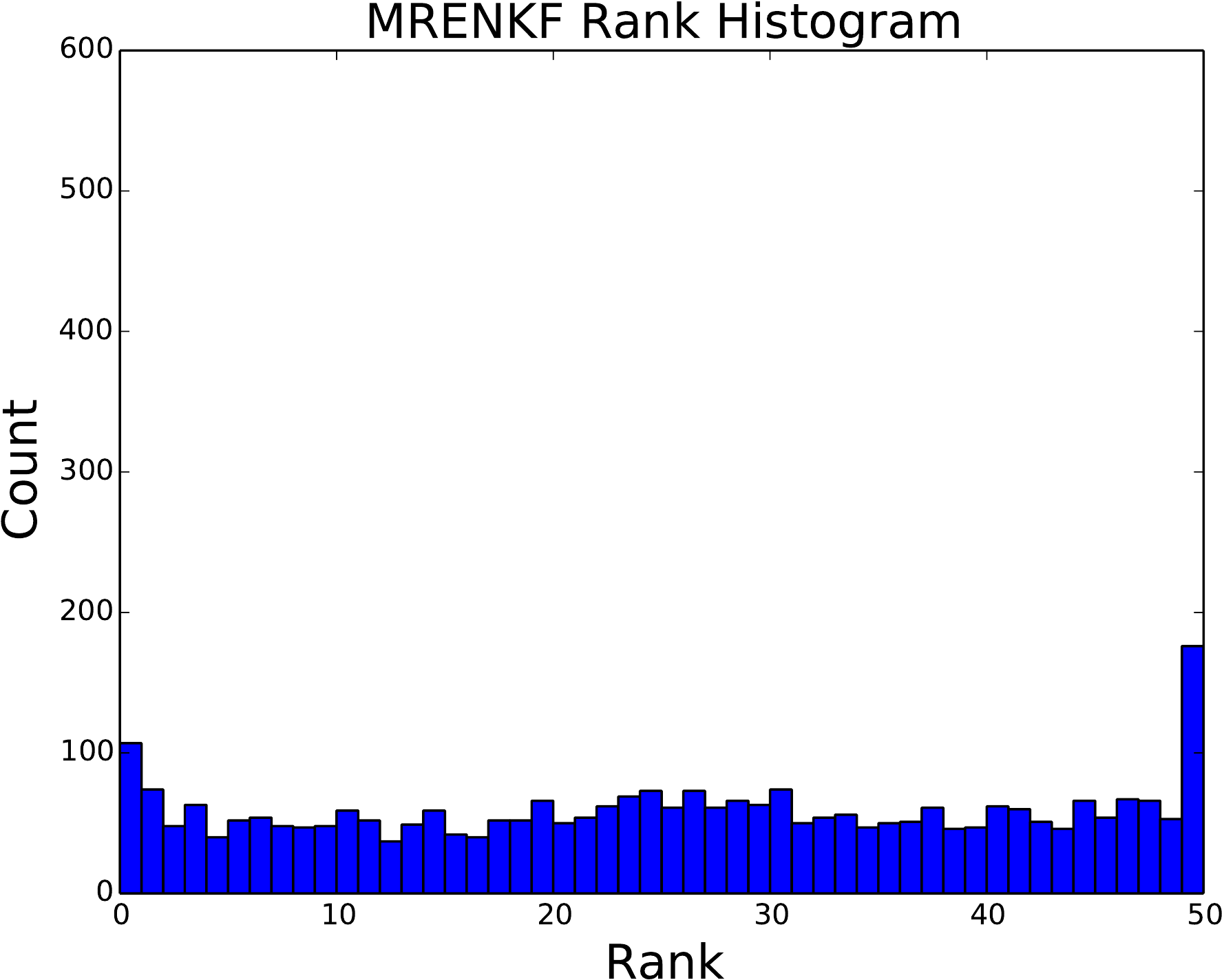}
      \end{array}$
  \end{center}
  \caption{Rank histograms are presented for the K--S assimilation
    investigation using the EnKF (left) and MrEnKF (right). These were
    generated by ranking the ensemble forecasts at $50$ equally spaced spatial
    points every $10^{th}$ time step, in order to reduce the effects of
    spatio-temporal correlations, and observing the rank that the true
    solution occurred in. For both the EnKF and MrEnKF these ranks where
    binned and the above histograms were formed. If the ensemble spread
    accurately represents the distribution of where the true solution is
    expected to fall given the observations then the rank histogram should be
    approximately uniform. In the above images we see that the EnKF rank
    histogram is very U-shaped, indicating an ensemble forecast that has
    collapsed. The rank histogram for the MrEnKF is much more uniform and
    therefore indicates a more accurate representation of the true solution's
    distribution given the observations.}
  \label{fig:etkfVSmretkf_rankhistogram}
\end{figure}

\begin{figure}[h]
  \begin{center}
      \includegraphics[scale=0.5]{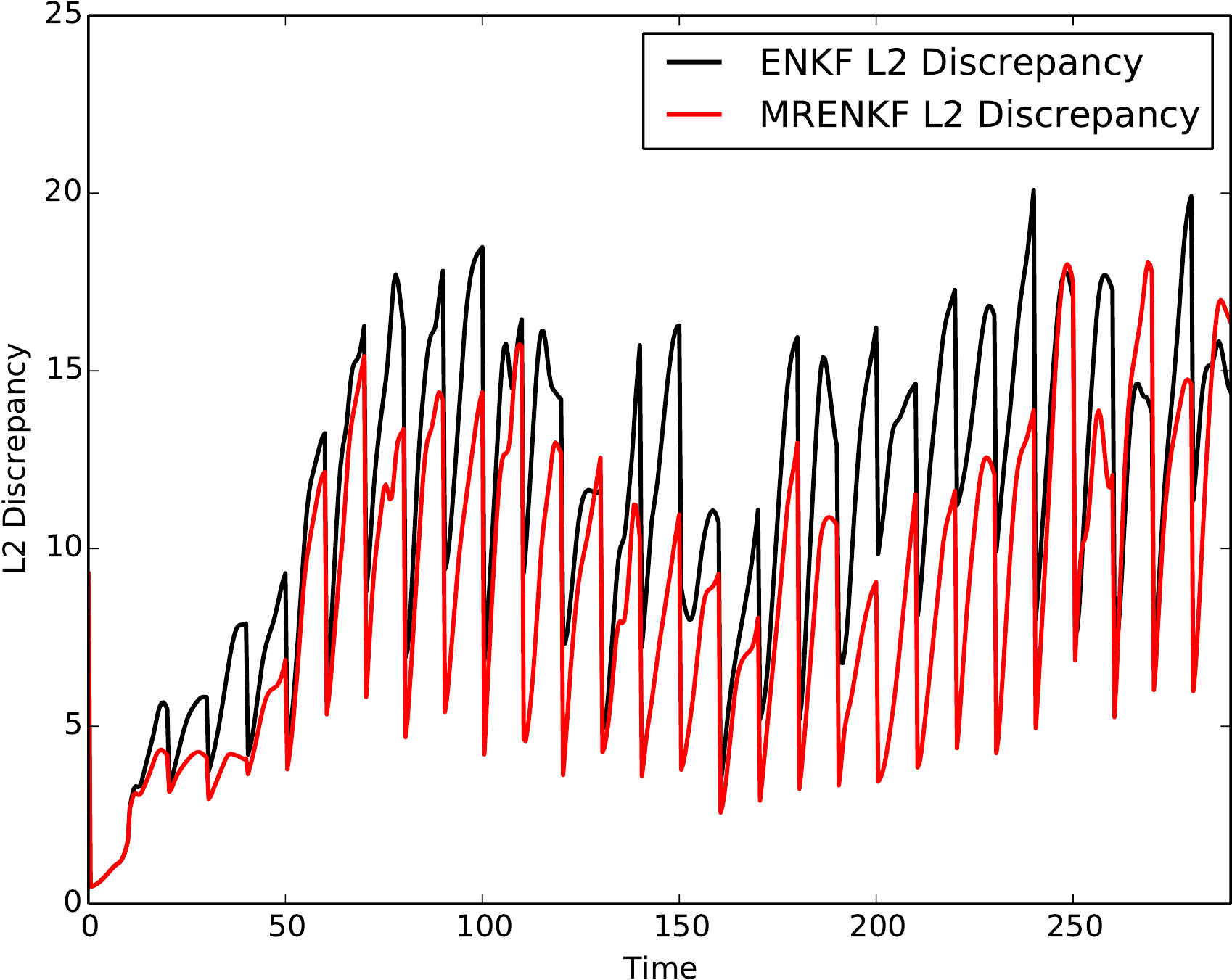}
  \end{center}
  \caption{Here we show the $L^2$ discrepancy between the true solution to the
    K--S equation and the mean of the ensemble forecast for the EnKF scheme
    (black) and MrEnKF scheme (red). We see that the MrEnKF scheme usually
    results in a lower forecast discrepancy. This fact coupled with the fact
    that the rank histogram is much more uniform for the MrEnKF, Figure
    \ref{fig:etkfVSmretkf_rankhistogram}, shows that the MrEnKF offers many
    advantages to the scale independent EnKF in this scenario. The discrepancy
    was calculated by computing the mean of the forecast at each time step and
    computing the $L^2$ difference between the mean and the true solution at
    every time step.}
  \label{fig:etkfVSmretkf_L2discrepancy}
\end{figure}

The EnKF rank histogram has a very U-shaped distribution indicating that the
true solution is most often found in the tails of the forecast distribution
(left plot in Figure~\ref{fig:etkfVSmretkf_rankhistogram}). On the other hand,
the MrEnKF rank histogram is close to uniform (right plot in
Figure~\ref{fig:etkfVSmretkf_rankhistogram}) indicating that the forecast
ensemble spread is a reliable representation of the true variation in the
solution given by the observed data.
For a secondary measure of performance evaluation we compute the $L^{2}$-norm
of the difference between the reference solution and the assimilation ensemble
mean, which we will refer to as the $L^2$-discrepancy for the mean.  In this
measure, we again see that the MrEnKF outperforms the standard EnKF
(Figure~\ref{fig:etkfVSmretkf_L2discrepancy}). The majority of the time, the
MrEnKF discrepancy is smaller than the EnKF discrepancy, which indicates that
MrEnKF provides a more accurate forecast. Both methods exhibit oscillations
in the $L^2$-discrepancy, where we see a sudden decrease in discrepancy where
we perform assimilation, and then a significant increase, indicating that the
model is deviating from the true solution. The deviation is due to the chaotic
dynamics within the K--S equation, which is reflected in the tendency of the
ensemble to spread out and away from the analysis after some time.

\subsection{Solar Photosphere}
\label{sec:ADAPT}

We apply the MrEnKF to a problem in solar weather using the Air Force Data
Assimilative Photospheric Flux Transport Model (ADAPT)
\cite{arge2013modeling,arge2010air,arge2011improving,hickmann2015data}.  In
ADAPT the magnetic flux is propagated across the Sun's surface using the
combined effects of differential rotation, meridional flow, and super granular
diffusion \cite{hickmann2015data,worden2000evolving}. ADAPT does well at
accurately transporting flux that is already present in the model ensemble's
forecast. However, solar physicists are interested in the tracking of emergent
coherent regions of magnetic flux of the same sign. These large clumps of
magnetic flux are known as \emph{active regions} and are primary drivers of
large space weather events such as Coronal Mass Ejections (CMEs)
\cite{antiochos1999model,falconer2002correlation,glover2000onset,munoz2010double}.
The underlying dynamic ADAPT model has no mechanism in place to generate these
active regions since the physics of their appearance is still not well
understood. Therefore, the ADAPT simulation has significant model error at the
scale of active regions. 

\begin{figure}[h]
  \begin{center}
    $\begin{array}{cc}
      \includegraphics[scale=0.4]{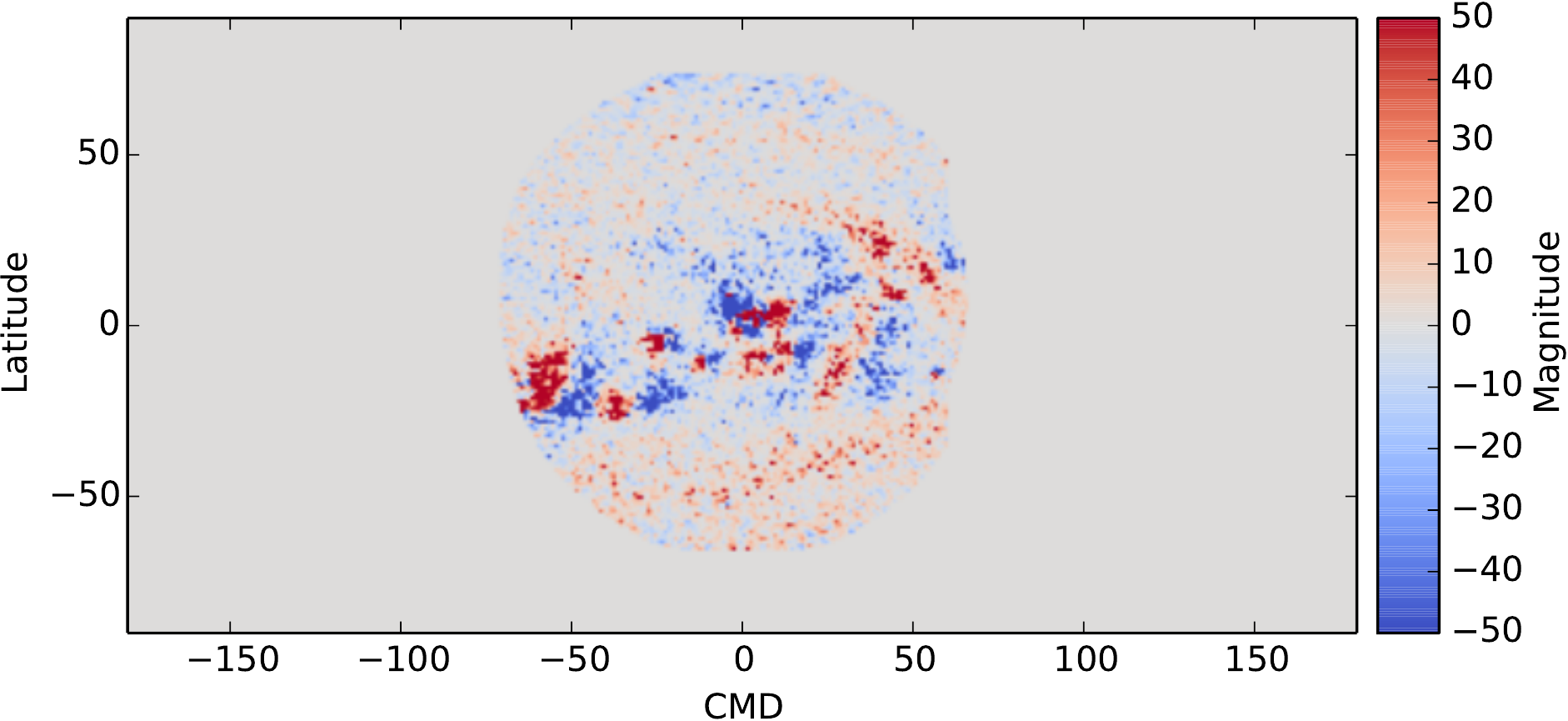} &
      \includegraphics[scale=0.4]{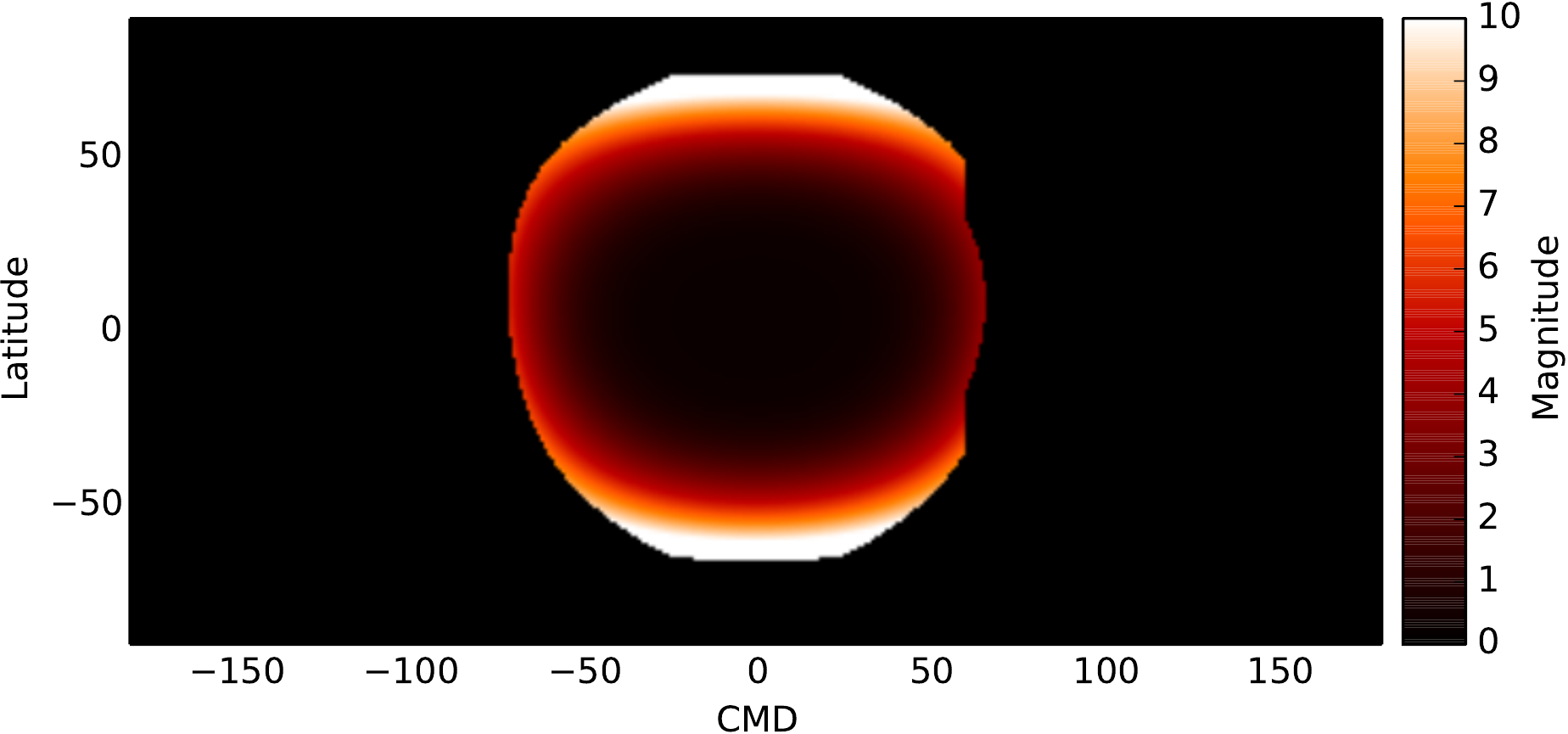} 
      \end{array}$
  \end{center}
  \caption{(LEFT) An example observation from the SOLIS vector
    spectro-magnetograph (VSM). The SOLIS-VSM provides observations of the
    radial magnetic flux from between $\pm 60^{\circ}$ latitude and $\pm
    90^{\circ}$ longitude, centered on the central meridian. The cadence of
    these observations is roughly one image every $24$ hours. These
    observations are then assimilated into ADAPT's evolving global map of the
    photospheric magnetic flux. Active regions are pictured as large coherent
    areas of magnetic flux with the same sign and are primary drivers of large
    scale solar weather events. As the sun rotates these active regions emerge
    on the \emph{east limb} (left side of the observation region) and then
    continue to be transported across the observation region. Since the radial
    magnetic flux is observed the curvature of the Sun near the boundaries of
    the observation region causes a larger observation error near the east
    limb. However, observations of active regions emerging on the east limb
    are trusted by solar weather experts and therefore should be assigned a
    high confidence during assimilation. (RIGHT) Pixel-by-pixel standard
    deviation associated with the SOLIS-VSM observation region (a standard
    deviation of $0$ corresponds to unobserved regions of the solar
    surface). The standard deviation in the center of the observation region
    are small but non-zero. Near the boundaries of the observation region the
    radial magnetic field has much higher observation error due to the
    curvature of the Sun's surface. This causes traditional EnKF methods to
    discard observations on the boundaries of the observation region. However,
    active regions near the boundaries are trusted a great deal by solar
    physicists and thus the MrEnKF serves an important purpose in assigning
    more trust to the observation of these large scale features.}
  \label{fig:solis_observation}
\end{figure}

ADAPT is updated using observations of the photospheric flux on the Earth side
of the Sun with an expert informed model of observation error, see Figure
\ref{fig:solis_observation}.
The model of observational noise used by solar physicists assumes
uncorrelated, zero mean, Gaussian measurement noise at each pixel with a
standard deviation that grows as the observation boundary is approached
\cite{henney2006solis,henney2007solis}. Observation error is greater near the
edges of the visible region of the Sun (right plot in
Figure~\ref{fig:solis_observation}), since only the radial component of the
magnetic flux is observed and the accuracy of this observation is proportional
to how aligned the observatory is with the direction of the radial component
at a point on the solar surface. Inevitably the edges of the visible portion
of the photosphere have radial directions not aligned with Earth and thus the
\emph{limbs} of the observation region are associated with greater error
\cite{henney2006solis,henney2007solis}.

This model of measurement error does not have any scale dependence. However,
strong fluctuations in the magnetic flux of the photosphere appear clumped
together in large \emph{active regions} \cite{munoz2010double}. When an active
region is observed, i.e. a large coherent region of magnetic flux all of the
same sign, solar physicists trust the observation and want to see that region
represented in the ensemble. This is especially important if the active region
is newly emerged into the observed region of the photosphere. Any data
assimilative algorithm used for magnetic flux forecasting should therefore
assign more trust to large coherent regions of magnetic flux with the same
sign and thus insert the observed active region into the analysis
ensemble. This can not be done by the standard EnKF with uncorrelated
pixel-by-pixel Gaussian measurement noise since the scale of an observed
feature does not effect its confidence during assimilation. A scale dependent
ensemble Kalman filter can make progress toward resolving this data
assimilation problem.

Unlike our previous example with the Kuramoto-Sivashinsky equation, in which
scale dependent observation noise necessitated the use of the MrEnKF, our
solar weather example has scale dependent \emph{model error}. In the ADAPT
model of photospheric flux dynamics there is no mechanism for the creation of
new active regions since these are caused by physical processes occurring
below the observable surface of the Sun. From the perspective of the solar
observation instruments the pixel-by-pixel representation of observation error
is sensible. However, since the model has no way to insert large scale active
regions, and these can arise over a short time span relative to the cadence of
solar observations, the model will necessarily be diverged from observations
of active regions that have appeared since the last observation. This effect
is especially apparent when active regions appear on the \emph{east limb}
(left side, Figure~\ref{fig:solis_observation}) of the observation region
since these solar regions have just emerged from the far side of the Sun and
have not been observed for a long span of time. The implication is that the
MrEnKF can serve a useful purpose in preparing solar weather forecasts by
systematically assigning a high confidence to observed structures at large
scales and a low confidence to observed structures occurring at small scales
within the ADAPT ensemble. This weighting is accomplished by reducing the
assumed observation error and increasing the ensemble inflation factor at
large scales.

\begin{figure}[h]
  \begin{center}
    \includegraphics[scale=0.42]{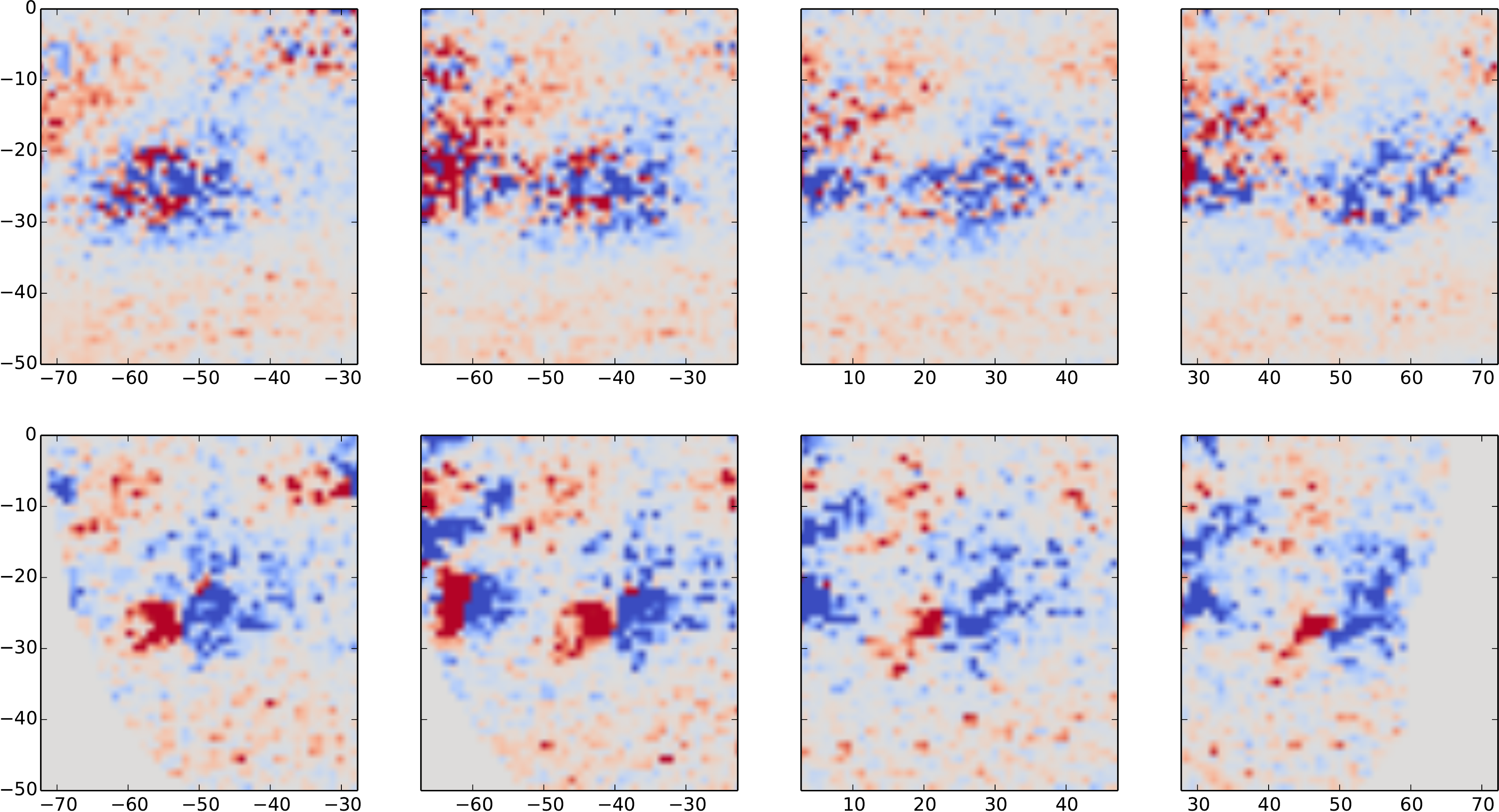}
  \end{center}
  \caption{Here we illustrate the effect of using a scale independent EnKF
    assimilation scheme on a small active region. The top row corresponds to
    the mean of the EnKF forecast while the bottom row corresponds to the
    SOLIS-VSM observations. Red and blue regions represent positive and
    negative polarity regions respectively, the active region in question
    appears as a horizontal red and blue pair centered in each frame of the
    observations. The $x$ and $y$ axis of each frame represent the location on
    the solar surface in latitude-longitude using \emph{Central Meridian
      Distance} (CMD) for the longitudinal coordinates. Observation times,
    from left to right, are the $18^{th}$ of November 2003 at 18:35, the
    $19^{th}$ of November 2003 at 18:11, the $24^{th}$ of November 2003 at
    17:05, and the $26^{th}$ of November 2003 at 17:49. We can see that the
    EnKF disperses the active region when it first appears on the East
    limb. This dissipation of the active region then continues as the active
    region is tracked across the observation region until the active region is
    almost completely dissipated in the forecast before it exits the
    observation region. This error greatly reduces the utility of the EnKF
    solar forecast.}
  \label{fig:etkf_activeregion}
\end{figure} 
\begin{figure}[h]
  \begin{center}
    \includegraphics[scale=0.42]{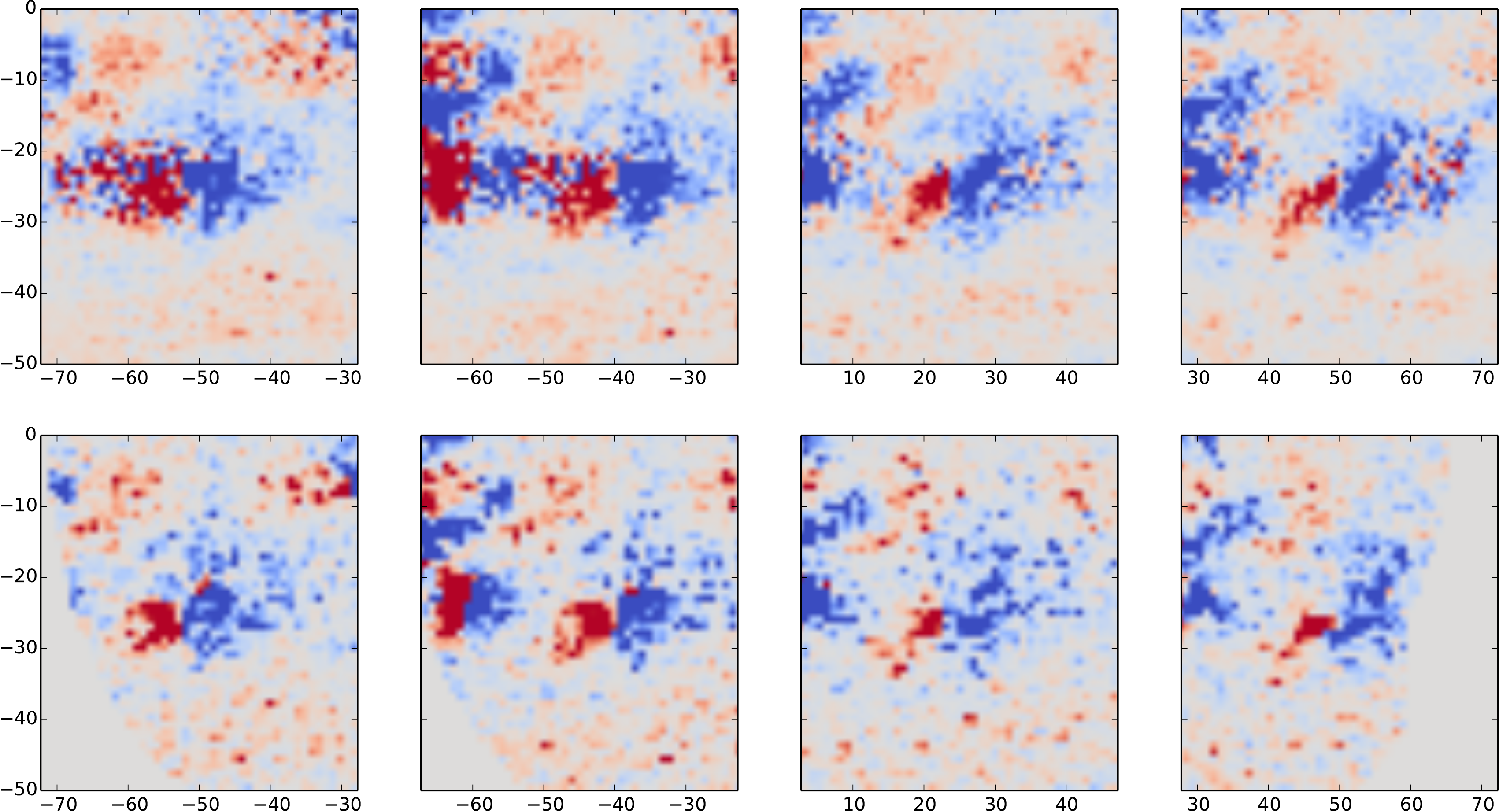}
  \end{center}
  \caption{Here we illustrate the effect of using our MrEnKF assimilation
    scheme on a small active region. The top row corresponds to the mean of
    the MrEnKF forecast while the bottom row corresponds to the SOLIS-VSM
    observations. Red and blue regions represent positive and negative
    polarity regions respectively, the active region in question appears as a
    horizontal red and blue pair centered in each frame of the
    observations. The $x$ and $y$ axis of each frame represent the location on
    the solar surface in latitude-longitude using \emph{Central Meridian
      Distance} (CMD) for the longitudinal coordinates. Observation times,
    from left to right, are the $18^{th}$ of November 2003 at 18:35, the
    $19^{th}$ of November 2003 at 18:11, the $24^{th}$ of November 2003 at
    17:05, and the $26^{th}$ of November 2003 at 17:49. We can see that
    the MrEnKF preserves the coherent structure of the active region since
    more confidence has been assigned to observations occurring at the scale
    of active regions. As the active region is then repeatedly observed during
    its journey across the observation region the MrEnKF refines the forecast
    structure of the region. Notice that the active region forecast is still
    well resolved as it exits the observation region on the West limb.}
  \label{fig:mretkf_activeregion}
\end{figure} 
Unlike our K--S example the ADAPT model state and observations exist in
$\mathbb{R}^2$ and therefore a two dimensional wavelet must be used. There
exist several ways to generalize one dimensional wavelets to higher
dimensions, in this work we use the tensor products
\cite{daubechies1992ten,mallat1989multiresolution,mallat1989theory} of
Daubechies 'db9' wavelets and only use a two level transformation on the
SOLIS-VSM observations. Figures~\ref{fig:etkf_activeregion} and
\ref{fig:mretkf_activeregion} compare the effect of using the EnKF and MrEnKF
to assimilate an active region emerging on the east limb of the solar
photosphere. The active region tracked in Figures~\ref{fig:etkf_activeregion}
and \ref{fig:mretkf_activeregion} emerged into the SOLIS-VSM observation
region on November $18^{th}$ 2003 at 18:35. Observation and assimilation of
this active region is shown as it passes across the observation region and
then exits on the West limb. We can see that the EnKF has difficulties
retaining scale dependent coherent features of the active region in
Figure~\ref{fig:etkf_activeregion}. The diffusion of the active region in the
EnKF algorithm becomes more pronounced as the active region is tracked across
the solar surface. As the active region exits the observation region and
crosses to the far side of the Sun the ensemble has almost completely diverged
from observations in the neighborhood of the active region. When we examine
the MrEnKF assimilation of the active region in
Figure~\ref{fig:mretkf_activeregion} we see that the method preserves a more
coherent structure of the active region. Moreover, as the active region is
tracked across the observation region it becomes more resolved. This is due to
the MrEnKF's ability to assign greater confidence to observed features at the
scale of active regions. Since, once the active region exits on the West limb,
it will not be observed again until it traverses the far side of the Sun it is
paramount to have a good estimation of the size and intensity of the active
region before it exits the observation region.

\section{Discussion}

We have detailed a method for inserting scale dependent information into an
ensemble Kalman filter framework. Our method was demonstrated on a 1D
nonlinear partial differential equation with scale dependent observation noise
and on an example from solar weather forecasting in which the model error, due
to missing physics, was highly scale dependent. The MrEnKF has the ability to
account for scale dependent variations in observation and model accuracy and
therefore tracked the evolving true solution of the Kuramoto--Sivashinsky
equation more accurately than the standard EnKF for small ensemble size. The
MrEnKF was also able to allow for scale dependent model deviations from
observations in our solar photosphere example.

The effect of scale dependent observations and model errors are common in many
scientific applications where forecasting is of interest. This problem is
usually handled in an \emph{ad hoc} way in practice by utilizing expert
opinions of the forecast's accuracy and manually adjusting observation error
accordingly on a point-by-point basis. By combining the EnKF with a
multi-scale wavelet analysis we have provided a general method to insert scale
dependent information, regarding model accuracy and observation accuracy, into
the assimilation scheme. Potential applications abound in areas of atmospheric
and oceanic forecasting in which models may be accurate for large scale,
non-turbulent effects, but highly unreliable at smaller scales.

In the Kuramoto-Sivashinsky application we demonstrated that, when the
observation error is scale dependent, and the ensemble size is small the
MrEnKF can be tuned to track the ground truth with less bias and more accuracy
than the EnKF. The rank histogram of the MrEnKF exhibited a more uniform
distribution than that of the EnKF, implying that the ensemble distribution
for the MrEnKF more accurately represented the observation probability
distribution. Multiresolution analysis combined with the EnKF was motivated,
for the authors, by the photospheric forecasting problem in solar weather. To
this end we have demonstrated the MrEnKF's effectiveness at preserving
coherent structures observed on the photosphere known as active regions. The
EnKF has difficulty capturing active region features in the analysis ensemble
since the underlying forecast model does not have physics to generate emerging
active regions. We plan to pursue a detailed study of active region
assimilation using the MrEnKF in forthcoming publications.

The MrEnKF method does need further development to be a readily
\emph{out-of-the-box} applicable tool. In particular the decision of which
wavelet basis to use and what level of wavelet transform to apply will be
explored in future work. In our applications several multiresolution levels
and wavelet types were experimented with. We note that the results did not
seem particularly sensitive to the choice of wavelet. Moreover, the level of
wavelet transform should be mostly informed by expert knowledge of the scale
at which observation accuracy shifts or model accuracy shifts. Currently we
solve this problem experimentally but plan to investigate adaptive methods in
future work.

Another research topic in fully developing the MrEnKF is to put forward a
reliable method to choose the ensemble inflation parameter at each scale. In
our examples we tuned the inflation at each scale through
experimentation. This is feasible if the number of multiresolution levels is
small but would become impractical for a high number of wavelet transform
levels. The methods pursued in
\cite{anderson2007adaptive,anderson2009spatially,li2009simultaneous,wang2003comparison,whitaker2012evaluating,ying2015adaptive}
involving an adaptive covariance inflation could circumvent this
difficulty. We intend to investigate these types of adaptive inflation schemes
in the context of our multiresolution EnKF in the future.

\section{Acknowledgments}

This research was primarily supported by NASA Living With a Star project
\#NNA13AB92I, ``Data Assimilation for the Integrated Global- Sun
Model''. Additional support was provided by the Air Force Office of Scientific
Research project R-3562-14-0, ``Incorporation of Solar Far-Side Active Region
Data within the Air Force Data Assimilative Photospheric Flux Transport
(ADAPT) Model''. The photospheric observations used in Figures
\ref{fig:solis_observation}, \ref{fig:etkf_activeregion}, and
\ref{fig:mretkf_activeregion} were provided by SOLIS-VSM. 


\bibliographystyle{plain}
\bibliography{MRwaveletDA_arXiv}

\end{document}